\documentclass[11pt]{amsart}
\usepackage{graphicx,amssymb,amsmath,amsthm}
\usepackage{enumerate}
\usepackage{comment,cite,color}
\usepackage{cite,color}
\usepackage{algorithmic}
\usepackage{mathrsfs}
\usepackage[ruled]{algorithm}
\usepackage{epsfig}
\usepackage{enumitem}

\newcommand{\innerp}[1]{\langle {#1} \rangle}

\newcommand{\abs}[1]{\lvert#1\rvert}

\newcommand{\argmin}[1]{\mathop{\rm argmin}\limits_{#1}}

\newcommand{\sign}{{\rm sign}}

\newcommand{\eproof}{\hfill\rule{2.2mm}{3.0mm}}

\newcommand{\Proof}{\noindent {\bf Proof.~~}}

\newcommand{\R}{{\mathbb R}}

\newcommand{\C}{{\mathbb C}}

\renewcommand{\eqref}[1]{(\ref{#1})}

\newcommand{\mhsp}{\hspace{2em}}

\newcommand{\supp}{{\rm supp}}

\newtheorem{prop}{Proposition}[section]

\newtheorem{definition}{Definition}[section]
\newtheorem{corollary}{Corollary}[section]
\newtheorem{theorem}{Theorem}[section]
\newtheorem{lemma}{Lemma}[section]
\newtheorem{remark}{Remark}[section]

\date{}

\begin{document}
\bibliographystyle{plain}
\title{Stable Signal Recovery from Phaseless Measurements }

\author{Bing Gao}
\address{LSEC, Inst.~Comp.~Math., Academy of
Mathematics and System Science,  Chinese Academy of Sciences, Beijing, 100091, China}
\email{gaobing@lsec.cc.ac.cn}

\author{Yang Wang}
\thanks{Yang Wang was supported in part by the
AFOSR grant FA9550-12-1-0455 and
NSF grant IIS-1302285.
       Zhiqiang Xu was supported  by NSFC grant 11171336, 11331012, 11021101 and National Basic Research
Program of China (973 Program 2010CB832702).}
\address{Department of Mathematics, HKUST}
\email{yangwang@ust.hk}

\author{Zhiqiang Xu}
\address{LSEC, Inst.~Comp.~Math., Academy of
Mathematics and System Science,  Chinese Academy of Sciences, Beijing, 100091, China}
\email{xuzq@lsec.cc.ac.cn}

 \maketitle{}

\begin{abstract}
The aim of this paper is to study the stability of the $\ell_1$ minimization for the compressive phase retrieval and  to extend
the instance-optimality in compressed sensing to the real phase retrieval setting. We first show that the
$m={\mathcal O}(k\log(N/k))$ measurements is enough to guarantee the
$\ell_1$ minimization  to recover $k$-sparse signals stably provided the measurement matrix $A$ satisfies
    the  strong RIP  property. We second  investigate  the phaseless instance-optimality with presenting a null space property of
    the measurement matrix $A$ under which there exists a decoder $\Delta$ so that the  phaseless instance-optimality holds.
We use the result to study the phaseless instance-optimality for the $\ell_1$ norm. The results build a parallel for compressive
 phase retrieval with the classical compressive sensing.
\end{abstract}

\section{Introduction}
In this paper we consider the phase retrieval for sparse signals with noisy measurements, which arises in many different applications. Assume that
 $$
 b_j:=\abs{\left<a_j,x_0\right>}+e_j, \quad j=1,\ldots,m
 $$
 where $x_0\in \R^N$, $a_j\in \R^N$ and $e_j\in \R$ is the noise.  Our goal is to recover $x_0$ up to a unimodular scaling constant from $b:=(b_1,\ldots,b_m)^\top$
 with the assumption of $x_0$ being approximately $k$-sparse. This problem is referred to as the {\em compressive phase retrieval problem} \cite{cpr07}.

The paper attempts to address two problems.  Firstly we consider the stability of $\ell_1$ minimization for the compressive phase retrieval problem where the signal $x_0$ is approximately $k$-sparse, which is the $\ell_1$ minimization problem defined as follows:
\begin{align}\label{opt:main}
     \min \|x\|_1\quad    \mbox{subject to} \quad \bigl\||Ax|-|Ax_0|\bigr\|_2\leq\epsilon,
\end{align}
where $A:=[a_1,\ldots,a_m]^\top$ and $\abs{Ax_0}:=[\abs{\left<a_1,x_0\right>},\ldots,\abs{\left<a_m,x_0\right>}]^\top$. Secondly we investigate the instance-optimality in the phase retrieval setting.

Note that in the classical compressive sensing setting  the stable recovering a $k$-sparse signal $x_0\in\C^N$ can be done using $m={\mathcal O}(k\log(N/k))$  measurements for several classes of measurement matrices $A$. A natural question is whether stable compressive phase retrieval can also be
  attained with $m={\mathcal O}(k\log(N/k))$ measurements. This is indeed proved to be the case in \cite{EM14} if $x_0\in\R^N$ and $A$ is
  a random real Gaussian matrix. One drawback is that the recovery method presented in \cite{EM14} is computationally unfeasible.
  In \cite{IVW14} a two-stage algorithm for compressive phase retrieval is proposed, which allows for very fast recovery of a sparse signal
  if the matrix $A$ can be written as a product of a random matrix and another matrix (such as a random matrix) that allows for efficient
   phase retrieval. The authors also proved that stable compressive phase retrieval can be achieved with  $m={\mathcal O}(k\log(N/k))$
    measurements also for complex signals $x_0$. In this paper, we shall show that the $\ell_1$ minimization scheme given in (\ref{opt:main})
    will recover a $k$-sparse signal stably from $m={\mathcal O}(k\log(N/k))$ measurements, provided that the measurement matrix $A$ satisfies
    the so-called strong RIP (S-RIP) property. This establishes an important parallel for compressive phase retrieval with the classical
    compressive sensing. Note that in \cite{WaXu14} such a parallel in terms of the null space property was already established.

The notion of {\em instance optimality} was first introduced in \cite{CDR09}. Given a norm $\|\cdot\|_X$ such as the $\ell_1$-norm and $x\in\R^N$, the best $k$-term approximation error is defined as
\[
\sigma_k(x)_X\,\,:=\,\, \min_{z\in \Sigma_k}\|x-z\|_X,
\]
where
$$
\Sigma_k:=\{x\in \R^N: \|x\|_0\leq k\}.
$$
We use $\Delta: \R^m\mapsto \R^N$ to denote a decoder for reconstructing $x$.  We say the pair $(A,\Delta)$ is {\em instance optimal of  order $k$ with constant $C_0$} if
\begin{equation}\label{eq:inst}
\|x-\Delta(Ax)\|_X\leq C_0\sigma_k(x)_X
\end{equation}
holds for all $x\in \R^N$. In extending it to phase retrieval, our decoder will have the input $b=\abs{Ax}$.  A pair $(A,\Delta)$ is said to be {\em phaseless instance optimal of order $k$ with constant $C_0$} if
\begin{equation}\label{eq:pinst}
 \min\Bigl\{\|x-\Delta(\abs{Ax})\|_X, \|x+\Delta(\abs{Ax})\|_X\Bigr\}\leq C_0\sigma_k(x)_X
\end{equation}
holds for all $x\in \R^N$. We are interested in the following problem : {\em Given $\|\cdot\|_X$ and $k<N$, what is the minimal value of $m$ for which there exists $(A,\Delta)$ so that (\ref{eq:pinst}) holds?}

The null space $\mathcal{N}(A):=\{x\in \R^N:Ax=0\}$ of $A$ plays an important role in the analysis of the original instance optimality (\ref{eq:inst}) (see \cite{CDR09}). Here we present a null space property for $\mathcal{N}(A)$, which is necessary and sufficient, for which there exists a  decoder $\Delta$ so that (\ref{eq:pinst}) holds. We apply the result to investigate the instance optimality where $X$ is the $\ell_1$ norm.
Set
\[
\Delta_1(\abs{Ax}):=\argmin{z\in \R^N}\Bigl\{\|z\|_1: \abs{Ax}=\abs{Az}\Bigr\}.
\]
We show that the pair $(A,\Delta_1)$ satisfies the (\ref{eq:pinst}) with $X=\ell_1$-norm provided $A$ satisfies the strong
RIP property (see Definition \ref{de:srip}). As shown in \cite{VX14}, the Gaussian random matrix $A\in \R^{m\times N}$ satisfies the strong RIP of order $k$ for
$m={\mathcal O}(k\log (N/k)$. Hence $m={\mathcal O}(k\log(N/k))$ measurements suffice to ensure the phaseless instance optimality (\ref{eq:pinst}) for  $\ell_1$-norm exactly as with the traditional
instance optimality (\ref{eq:inst}).

\section{Auxiliary Results}
\setcounter{equation}{0}

In this section we provide some auxiliary results that will be used in later sections.
For $ x\in\R^N $ we use $\|x\|_p:=\|x\|_{\ell_p}$ to denote the $p$-norm of $x$ for $0<p \leq \infty$.
The measurement matrix is given by $A:=[a_1,\ldots,a_m]^T \in\mathbb{R}^{m\times N}$ as before. Given an index set $I\subset \{1,\ldots,m\}$  we shall use  $A_I$ to denote the sub-matrix of $A$ where only rows with indices in $I$ are kept, i.e.,
\[
A_I:=[a_j:j\in I]^\top.
\]
The matrix $A$ satisfies the {\em Restricted Isometry Property (RIP) of order $k$} if there exists a constant  $\delta_k\in [0,1)$ such that for all $k$-sparse vectors $z\in\mathbb{R}^N$ we have
$$
(1-\delta_k)\|z\|_2^2\leq\|Az\|_2^2\leq(1+\delta_k)\|z\|_2^2.
$$
It was shown in \cite{CZ14} that one can use $\ell_1$ minimization to recover $k$-sparse signals provided that $A$ satisfies the RIP of order $t\cdot k$ and $\delta_{t\cdot k}<\sqrt{1-\frac{1}{t}}$ where $t>1$.

To investigate compressive phase retrieval, a stronger notion of RIP is given in \cite{VX14}:

\begin{definition}(S-RIP)\label{de:srip}:
We say the matrix $A=[a_1,\cdots,a_m]^\top \in\mathbb{R}^{m\times N}$  has the {\em Strong Restricted Isometry Property} of order k with bounds $\theta_-,\ \theta_+\in (0, 2)$ if
\begin{equation} \label{eq:SRIP}
\theta_-\|x\|_2^2\leq \min_{ I\subseteq[m], |I|\geq m/2}\|A_Ix\|_2^2\leq\max_{I\subseteq[m],|I|\geq m/2}
\|A_Ix\|_2^2\leq\theta_+\|x\|_2^2
\end{equation}
holds for all k-sparse signals $x\in\mathbb{R}^N$, where $[m]:=\{1,\ldots,m\}$. We say $A$ has the {\em Strong Lower Restricted Isometry Property} of order k with bound $\theta_-$ if the lower bound in (\ref{eq:SRIP}) holds. Similarly we say $A$ has the {\em Strong Upper Restricted Isometry Property} of order k with bound $\theta_+$ if the upper bound in (\ref{eq:SRIP}) holds.
\end{definition}

The authors of \cite{VX14} proved that Gaussian matrices with $m=\mathcal{O}(tk\log(N/k))$ satisfy S-RIP of order $tk$ with high probability.

\begin{theorem}[\cite{VX14}]\label{gaussiansrip}
Suppose that $t>1$ and $ A=(a_{ij})\in \mathbb{R}^{m\times N} $ is a random Gaussian matrix with $m=\mathcal{O}(tk\log(N/k))$
and $a_{ij}\sim {\mathcal N}(0,\frac{1}{\sqrt{m}})$. Then there exist $\theta_-, \theta_+ \in (0,2)$ such that with probability $1-\exp(-cm/2)$ the matrix $ A $ satisfies the S-RIP of order $ tk $ with constants  $\theta_-$ and   $\theta_+$, where $c>0$ is an absolute constant and $\theta_-$,  $\theta_+$ are independent of $t$.
\end{theorem}

The following is a very useful lemma for this study.

\begin{lemma}  \label{mainlemma}
 Suppose that $ x_0\in\mathbb{R}^N$ and  $ \rho\geq 0$. Suppose that $ A\in\mathbb{R}^{m\times N}$ is a measurement matrix satisfying the restricted isometry property with $ \delta_{tk}\leq\sqrt{\frac{t-1}{t}} $ for some $ t>1 $. Then for any
  \[
   \hat{x}\in \Bigl\{x\in \R^N : \|x\|_1\leq \|x_0\|_1+\rho, \, \|Ax-Ax_0\|_2\leq\epsilon \Bigr\}
  \]
 we have
 \begin{equation*}
   \|\hat{x}-x_0\|_2\leq c_1\epsilon+c_2\frac{2\sigma_k(x_0)_1}{\sqrt{k}}+c_2\cdot\frac{\rho}{\sqrt{k}},
  \end{equation*}
  where $ c_1=\frac{\sqrt{2(1+\delta)}}{1-\sqrt{t/(t-1)}\delta} $, $  c_2=\frac{\sqrt{2}\delta+\sqrt{(\sqrt{t(t-1)}-\delta t)\delta}}{\sqrt{t(t-1)}-\delta t}+1.$
\end{lemma}

\begin{remark}{\rm
We build the proof of Lemma \ref{mainlemma} following the ideas of Cai and Zhang \cite{CZ14}.The full proof is given in Appendix for completeness. It is well-known that an effective method to recover approximately-sparse signals $x_0$ in the traditional compressive sensing is to solve
\begin{equation}\label{eq:xjing}
    x^\#:=\argmin x\{\|x\|_1 : \|Ax-Ax_0\|_2\leq \epsilon\}.
\end{equation}
The definition of $x^\#$ shows  that
 $$
 \|x^\#\|_1\leq \|x_0\|_1, \|Ax^\#-Ax_0\|_2\leq \epsilon,
 $$
 which implies that
 $$
 \|x^\#-x_0\|_2\leq C_1\epsilon+C_2\frac{\sigma_k(x_0)_1}{\sqrt{k}},
 $$
 provided that $A$ satisfies the RIP condition with $\delta_{tk}\leq \sqrt{1-1/t}$ for $t>1$ (see \cite{CZ14}). However, in practice one prefers to design fast algorithms to find an approximation solution of (\ref{eq:xjing}), say $\hat{x}$. Thus it is possible to have
 $\|\hat{x}\|_1> \|x_0\|_1$. Lemma \ref{mainlemma} gives an estimation of $\|\hat{x}-x_0\|_2$ for the case where
 $\|\hat{x}\|_1\leq \|x_0\|_1+\rho$.
}
\end{remark}

\section{Stable Recovery of Real Phase Retrieval Problem}
\setcounter{equation}{0}

\subsection{Stability Results}
The following lemma shows that the map $\phi_A(x):=\abs{Ax}$ is stable on $\Sigma_k$ modulo a unimodular constant provided $A$ satisfies strong lower RIP of order $2k$. Define the equivalent relation $\sim$ on $\R^N$ and $\C^N$ by the following: for any $x \sim y$ iff $x= cy$ for some unimodular scalar $c$, where $x, y$ are in $\R^N$ or $\C^N$. For any subset $Y$ of $\R^N$ or $\C^N$   the notation $Y/\sim$ denotes the equivalent classes of elements in $Y$ under the equivalence. Note that there is a natural metric $D_\sim$ on $\C^N/\sim$ given by
$$
     D_\sim(x, y) = \min_{|c|=1} \|x-cy\|.
$$
Our primary focus in this paper will be on $\R^N$, and in this case
$D_\sim(x,y) = \min\{\|x-y\|_2, \|x+y\|_2\}$.

\begin{lemma}\label{th:stab}
   Let $A\in \R^{m\times N}$ satisfy the strong lower RIP of order $2k$ with constant $\theta_-$. Then for any $x, y \in \Sigma_k$ we have
$$
   \|\abs{Ax}-\abs{Ay}\|_2^2 \geq \theta_- \min (\|x-y\|_2^2, \|x+y\|_2^2).
$$
\end{lemma}
\Proof
For any $x,y\in \Sigma_k $ we divide $\{1,\ldots,m\}$ into two groups:
  $T=\{j:~\sign(\innerp{a_j,x})=\sign(\innerp{a_j,y})\}$ and $T^c=\{j:~\sign(\innerp{a_j,x})=-\sign(\innerp{a_j,y})\}$. Clearly one of $T$ and $T^c$ will have cardinality at least $m/2$. Without loss of generality we assume that $T$ has cardinality no less than $m/2$. Then
\begin{eqnarray*}
\|\abs{Ax}-\abs{Ay}\|^2_2 &=& \|A_Tx-A_T y\|_2^2 + \|A_{T^c}x+A_{T^c}y\|_2^2\\
  &\geq& \|A_Tx-A_T y\|_2^2 \\
  &\geq& \theta_- \|x-y\|_2^2 \\
  &\geq& \theta_- \min (\|x-y\|_2^2, \|x+y\|_2^2).
\end{eqnarray*}

\begin{remark}
Note that the combination of Lemma \ref{th:stab} and Theorem \ref{gaussiansrip} shows that for an $m\times N$ Gaussian matrix $A$ with $m=O(k\log (N/k))$ one can guarantee the map $\phi_A(x):=|Ax|$ is stable on $\Sigma_k/\sim$.
\end{remark}

\subsection{The Main Theorem}
In this part, we will consider how many measurements are needed for the stable sparse phase retrieval by
$\ell_1$-minimization via solving the following model:
\begin{equation} \label{opt:section2}
    \min\|x\|_1\mhsp \mbox{subject to} \mhsp \||Ax|-|Ax_0|\|_2^2\leq\epsilon^2,
\end{equation}
where $A$ is our measurement matrix and $x_0\in\mathbb{R}^N$ is a signal we wish to recover. The next theorem tells under what conditions the solution to (\ref{opt:section2}) is stable.

\begin{theorem}\label{maintherom}
Assume that $A\in \mathbb{R}^{m\times N}$ satisfies the S-RIP of order $tk$  with bounds $\theta_-, \theta_+ \in(0,2)$ such that
$$
       t\geq\max\{\frac{1}{2\theta_--\theta_-^2},\frac{1}{2\theta_+-\theta_+^2}\}.
$$
Then any solution $\hat x$ for (\ref{opt:section2}) satisfies
$$
   \min\{\|\hat{x}-x_0\|_2,\|\hat{x}+x_0\|_2\}\leq c_1\epsilon+c_2\frac{2\sigma_k(x_0)_1}{\sqrt{k}},
$$
where $c_1$ and $ c_2$ are constants defined in Lemma \ref{mainlemma}.
\end{theorem}
\Proof   Clearly any $\hat x\in \R^N$ satisfying (\ref{opt:section2}) must have
\begin{equation}\label{option1}
    \|\hat{x}\|_1\leq\|x_0\|_1
\end{equation}
and
\begin{equation}\label{eq:option2}
   \||A\hat{x}|-|Ax_0|\|_2^2\leq\epsilon^2.
\end{equation}
Now the indices set $\{1, 2, \dots, m\}$ is divided into two subsets
\begin{align*}
  T & =\{j:~ \sign(\innerp{a_j,\hat{x}})=\sign(\innerp{a_j,x_0})\} \\
  T^c & =\{j:~ \sign(\innerp{a_j,\hat{x}})=-\sign(\innerp{a_j,x_0})\}.
\end{align*}
Then (\ref{eq:option2}) implies that
\begin{equation}\label{eq:option3}
   \|A_T\hat{x}-A_Tx_0\|_2^2+\|A_{T^c}\hat{x}+A_{T^c}x_0\|_2^2\leq\epsilon^2.
\end{equation}
Here either $|T|\geq m/2$ or  $|T^c|\geq m/2$. Without loss of generality we assume that  $|T|\geq m/2$. We use the fact
\begin{equation}\label{eq:option4}
          \|A_T\hat{x}-A_Tx_0\|_2^2\leq\epsilon^2.
\end{equation}
From (\ref{option1}) and (\ref{eq:option4}) we obtain
  \begin{equation}\label{eq:set1}
  \hat{x}\in \{x\in \R^N: \|x\|_1\leq \|x_0\|_1, \|A_Tx-A_Tx_0\|_2\leq \epsilon\}.
  \end{equation}
Recall that $A$ satisfies S-RIP of order $tk$ and constants $\theta_-, \ \theta_+$. Here
\begin{equation}\label{eq:ttheta}
    t\geq\max \{\frac{1}{2\theta_--\theta_-^2},\frac{1}{2\theta_+-\theta_+^2}\}>1.
\end{equation}
The definition of S-RIP implies that  $A_T$ satisfies the RIP of order $tk$ in which
   \begin{equation}\label{eq:set2}
   \delta_{tk} \leq\max\{1-\theta_-,\ \theta_+-1\}\leq \sqrt{\frac{t-1}{t}}
   \end{equation}
   where the second inequality follows from (\ref{eq:ttheta}).
The combination of (\ref{eq:set1}), (\ref{eq:set2}) and Lemma \ref{mainlemma} now implies
$$
   \|\hat{x}-x_0\|_2\leq c_1\epsilon+c_2\frac{2\sigma_k(x_0)_1}{\sqrt{k}},
$$
where $c_1$ and $c_2$ are defined in Lemma \ref{mainlemma}. If $|T^c|\geq\frac{m}{2}$ we get the corresponding result
 \begin{equation*}
 \|\hat{x}+x_0\|_2\leq  c_1\epsilon+c_2\frac{2\sigma_k(x_0)_1}{\sqrt{k}}.
 \end{equation*}
The theorem is now proved.
\eproof

This theorem demonstrates that, if the measurement matrix have the S-RIP, real compressive phase retrieval problem can be solved stably by $\ell_1 $-minimization.

\section{Phase Retrieval and Best k-term Approximation}

\setcounter{equation}{0}

\subsection{Instance optimality from the linear measurements}
We introduce some definitions and results in \cite{CDR09}. Recall that for a given encoder matrix $A\in \R^{m\times N}$ and a decoder $\Delta:\R^m\mapsto \R^N$, the pair $(A,\Delta)$ is said to have instance optimality of  order $k$ with constant $C_0$ with respect to the norm $X$ if
\begin{equation}\label{eq:inst2}
    \|x-\Delta(Ax)\|_X\leq C_0\sigma_k(x)_X
\end{equation}
holds for all $x\in \R^N$. Set
   ${\mathcal N}(A):=\{\eta\in \R^N: A\eta=0\}$ to be the null space of $A$. The following theorem gives conditions under which the  (\ref{eq:inst2}) holds.

\begin{theorem}[\cite{CDR09}] \label{th:null}
   Let $A\in \R^{m\times N}$, $1 \leq k \leq N$ and $\|\cdot\|_X$ be a norm. Then a
sufficient condition for the existence of a decoder $\Delta$ satisfying (\ref{eq:inst2}) is
\begin{equation}\label{eq:null}
    \|\eta\|_X\leq \frac{C_0}{2}\sigma_{2k}(\eta)_X, \quad \forall \eta\in {\mathcal N}(A).
\end{equation}
A necessary condition for the existence of a decoder $\Delta$ satisfying (\ref{eq:inst2}) is
\begin{equation}\label{eq:necnull}
     \|\eta\|_X\leq C_0\sigma_{2k}(\eta)_X, \quad \forall \eta\in {\mathcal N}(A).
\end{equation}
\end{theorem}

For the norm $X=\ell_1$ it was established in \cite{CDR09} that instance optimality of order $k$ can indeed be achieved, e.g. for a Gaussian matrix $A$, with $m=O(k\log(N/k))$. The authors also considered more generally taking different norms on both sides of (\ref{eq:inst2}). Following \cite{CDR09}, we say the pair $(A,\Delta)$ has {\em $(p,q)$-instance optimality of  order $k$ with constant $C_0$} if
\begin{equation}\label{eq:pqinst2}
   \|x-\Delta(Ax)\|_{\ell_p}\leq C_0\sigma_k(x)_{\ell_q}/k^{1/p-1/q}, \quad \forall x\in \R^N,
\end{equation}
with $1\leq q\leq p \leq 2$. It was shown in \cite{CDR09} that the $(p,q)$-instance optimality of order $k$ can be achieved at the cost of having $m=\mathcal{O}(k(N/k)^{2-2/q})\log (N/k)$ measurements.

\subsection{Phaseless Instance Optimality }
     A natural question here is whether an analogous result of Theorem \ref{th:null} exists for phaseless instance optimality defined in (\ref{eq:pinst}). We answer the question by presenting such a result in the case of real phase retrieval.

\begin{theorem}\label{nsp_general}
    Let $A\in \R^{m\times N}$, $1 \leq k \leq N$ and $\|\cdot\|_X$ be a norm. Then a
sufficient condition for the existence of a decoder $\Delta$ satisfying the phaseless instance optimality (\ref{eq:pinst}) is: For any $ I\subseteq\{1,\ldots,m\}$ and $\eta_1\in\mathcal{N}(A_I)$, $\eta_2\in\mathcal{N}(A_{I^c})$ we have
\begin{equation}\label{eq:thsuff}
    \min\{\|\eta_1\|_X  ,  \|\eta_2\|_X \}\leq\frac{C_0}{4}\sigma_k(\eta_1-\eta_2)_X+\frac{C_0}{4}\sigma_k(\eta_1+\eta_2)_X.
\end{equation}
A necessary condition for the existence of a decoder $\Delta$ satisfying (\ref{eq:pinst}) is: For any $I\subseteq\{1,\ldots,m\}$ and $\eta_1\in\mathcal{N}(A_I)$, $\eta_2\in\mathcal{N}(A_{I^c})$ we have
\begin{equation}\label{eq:thnec}
     \min\{\|\eta_1\|_X  ,  \|\eta_2\|_X \}\leq\frac{C_0}{2}\sigma_k(\eta_1-\eta_2)_X+\frac{C_0}{2}\sigma_k(\eta_1+\eta_2)_X.
\end{equation}
\end{theorem}
\Proof
We first assume (\ref{eq:thsuff}) holds, and show that there exists a decoder $ \Delta $ satisfying the phaseless instance optimality (\ref{eq:pinst}). To this end, we define a decoder $ \Delta $ as follows:
$$
\Delta(|Ax_0|)=\mathop{\textup{argmin}}_{|Ax|=|Ax_0|}\sigma_k(x)_X.
$$
Suppose $ \hat{x}:=\Delta(|Ax_0|)$.  We have $|A\hat{x}|=|Ax_0|$ and $\sigma_k(\hat{x})_X\leq\sigma_k(x_0)_X$. Note that $\innerp{a_j,\hat{x}}=\pm \innerp{a_j,x_0}$. Let $I \subseteq\{1,\ldots,m\}$ be defined by
$$
  I=\Bigl\{j:~\innerp{a_j,\hat{x}} = \innerp{a_j,x_0}\Bigr\}.
$$
Then
$$
     A_I(x_0-\hat{x})=0,\quad   A_{I^c}(x_0+\hat{x})=0.
$$
Set
\begin{align*}
 \eta_1&:=x_0-\hat{x}\in\mathcal{N}(A_I),\\
 \eta_2&:=x_0+\hat{x}\in\mathcal{N}(A_{I^c}).
\end{align*}
A simple observation yields
\begin{equation}\label{eq:opt32}
       \sigma_k(\eta_1-\eta_2)_X=2\sigma_k(\hat{x})_X\leq2\sigma_k(x_0)_X , \quad
       \sigma_k(\eta_1+\eta_2)_X=2\sigma_k(x_0)_X.
\end{equation}
Then (\ref{eq:thsuff}) implies that
\begin{align*}
\min\{ \|\hat{x}-x_0\|_X, \|\hat{x}+x_0\|_X \}
           &= \min\{\|\eta_1\|_X  ,  \|\eta_2\|_X \} \\
           &\leq\frac{C_0}{4}\sigma_k(\eta_1-\eta_2)_X+\frac{C_0}{4}\sigma_k(\eta_1+\eta_2)_X \\
           &\leq C_0\sigma_k(x_0)_X.
\end{align*}
Here the last equality is obtained by (\ref{eq:opt32}). This proves the sufficient condition.

We next turn to the necessary condition.  Let $\Delta $ be a decoder for which the phaseless instance optimality (\ref{eq:pinst}) holds. Let $I\subseteq\{1,\ldots,m\}$. For any $\eta_1\in\mathcal{N}(A_I)$ and $\eta_2\in\mathcal{N}(A_{I^c})$ we have
\begin{equation} \label{opt3:3}
    |A(\eta_1+\eta_2)|=|A(\eta_1-\eta_2)|=|A(\eta_2-\eta_1)|.
\end{equation}
The instance optimality implies
\begin{equation}\label{eq:li1}
    \min\{\|\Delta(|A(\eta_1+\eta_2)|)+\eta_1+\eta_2\|_X,\|\Delta(|A(\eta_1+\eta_2)|)-(\eta_1+\eta_2)\|_X\}
              \leq C_0\sigma_k(\eta_1+\eta_2)_X.
 \end{equation}
 Without loss of generality we may assume that
$$
 \|\Delta(|A(\eta_1+\eta_2)|)+\eta_1+\eta_2\|_X\,\,\leq\,\,\|\Delta(|A(\eta_1+\eta_2)|)-(\eta_1+\eta_2)\|_X.
$$
Then (\ref{eq:li1}) implies that
\begin{equation}  \label{opt3:4}
   \|\Delta(|A(\eta_1+\eta_2)|)+\eta_1+\eta_2\|_X \leq C_0\sigma_k(\eta_1+\eta_2)_X.
\end{equation}
By (\ref{opt3:3}), we have
\begin{align}
       \|\Delta(|A(\eta_1+\eta_2)|)&+\eta_1+\eta_2\|_X\nonumber\\
             &=\|\Delta(|A(\eta_2-\eta_1)|)-(\eta_2-\eta_1)+2\eta_2\|_X \nonumber\\
             &\geq2\|\eta_2\|_X-\|\Delta(|A(\eta_2-\eta_1)|)-(\eta_2-\eta_1)\|_X.\label{opt3:5}
\end{align}
Combining (\ref{opt3:4}) and (\ref{opt3:5}) yields
\begin{equation} \label{opt3:6}
    2\|\eta_2\|_X\leq C_0\sigma_k(\eta_1+\eta_2)_X+\|\Delta(|A(\eta_2-\eta_1)|)-(\eta_2-\eta_1)\|_X.
\end{equation}
At the same time, (\ref{opt3:3}) also implies
\begin{align}
       \|\Delta(|A(\eta_1+\eta_2)|)&+\eta_1+\eta_2\|_X\nonumber\\
             &=\|\Delta(|A(\eta_2-\eta_1)|)+(\eta_2-\eta_1)+2\eta_1\|_X \nonumber\\
             &\geq2\|\eta_1\|_X-\|\Delta(|A(\eta_2-\eta_1)|)+(\eta_2-\eta_1)\|_X.\label{opt3:7}
\end{align}
Putting (\ref{opt3:4}) and (\ref{opt3:7}) together, we obtain
\begin{equation}  \label{opt3:8}
  2\|\eta_1\|_X\leq C_0\sigma_k(\eta_1+\eta_2)_X+\|\Delta(|A(\eta_2-\eta_1)|)+(\eta_2-\eta_1)\|_X.
\end{equation}
It follows from (\ref{opt3:6}) and (\ref{opt3:8}) that
\begin{align*}
\min\left\{\|\eta_1\|_X,\|\eta_2\|_X\right\}
   &\leq \frac{C_0}{2}\sigma_k(\eta_1+\eta_2)_X+\\
   &\frac{1}{2}\min\{
   \|\Delta(|A(\eta_2-\eta_1)|)-(\eta_2-\eta_1)\|_X,\|\Delta(|A(\eta_2-\eta_1)|)+(\eta_2-\eta_1)\|_X\}\\
   &\leq\frac{C_0}{2}\sigma_k(\eta_1+\eta_2)_X+\frac{C_0}{2}\sigma_k(\eta_1-\eta_2)_X.
\end{align*}
Here the last inequality is obtained by $ (A,\Delta) $ satisfying the instance optimality.
For the case where
\[
\|\Delta(|A(\eta_1+\eta_2)|)-(\eta_1+\eta_2)\|_X\,\,\leq\,\,\|\Delta(|A(\eta_1+\eta_2)|)+\eta_1+\eta_2\|_X,
\]
we obtain
\[
\min\{\|\eta_1\|_X,\|\eta_2\|_X\}\leq\frac{C_0}{2}\sigma_k(\eta_1+\eta_2)_X+\frac{C_0}{2}\sigma_k(\eta_1-\eta_2)_X
\]
by the similar way.
The theorem is now proved.
\eproof


We next present a null space property for phaseless instance optimality, which allows us to establish parallel results for sparse phase retrieval.

\begin{definition}{\rm
   We say a matrix $A \in \R^{m\times N}$ satisfies the {\em strong null space property (S-NSP) of order $ k $ with constant $ C $} if for any index set $I\subseteq\{1,\ldots,m\}$ with $\abs{I}\geq m/2$ and $ \eta \in {\mathcal N}(A_I)$ we have
$$
    \|\eta\|_X\leq C\cdot \sigma_{k}(\eta)_X.
$$
}
\end{definition}

\begin{theorem}\label{lemmaSNSP}
 Assume that a matrix $ A\in\mathbb{R}^{m\times N}$ has the strong null space property of order $ 2k $ with constant $ C_0/2 $. Then there must exist a decoder $\Delta$ having the phaseless instance optimality (\ref{eq:pinst}) with constant $ C_0 $. In particular, one such decoder is
$$
  \Delta(|Ax_0|)=\mathop{\textup{argmin}}_{|Ax|=|Ax_0|}\sigma_k(x)_X.
$$
\end{theorem}
\Proof
   Assume that $ I\subseteq\{1,\ldots,m\}$. For any $\eta_1\in\mathcal{N}(A_I)$ and $\eta_2\in\mathcal{N}(A_{I^c})$ we must have either $\|\eta_1\|_X\leq \frac{C_0}{2}\sigma_{2k}(\eta_1)_X$ or   $\|\eta_2\|_X\leq \frac{C_0}{2}\sigma_{2k}(\eta_2)_X$ by the strong null space property. If $\|\eta_1\|_X\leq \frac{C_0}{2}\sigma_{2k}(\eta_1)_X$ then
$$
  \|\eta_1\|_X\leq  \frac{C_0}{2}\sigma_{2k}(\eta_1)_X\leq
        \frac{C_0}{4}\sigma_k(\eta_1-\eta_2)_X+\frac{C_0}{4}\sigma_k(\eta_1+\eta_2)_X.
$$
Similarly if $\|\eta_2\|_X\leq \frac{C_0}{2}\sigma_{2k}(\eta_2)_X$ we will have
$$
   \|\eta_2\|_X\leq  \frac{C_0}{2}\sigma_{2k}(\eta_2)_X\leq
        \frac{C_0}{4}\sigma_k(\eta_1-\eta_2)_X+\frac{C_0}{4}\sigma_k(\eta_1+\eta_2)_X.
$$
It follows that
\begin{equation}\label{eq:1thsuff}
  \min\{\|\eta_1\|_X  ,  \|\eta_2\|_X \}
     \leq\frac{C_0}{4}\sigma_k(\eta_1-\eta_2)_X+\frac{C_0}{4}\sigma_k(\eta_1+\eta_2)_X.
\end{equation}
Theorem \ref{nsp_general} now implies that the required decoder $\Delta$ exists. Furthermore, by the proof of the sufficiency part of Theorem \ref{nsp_general},
  $$
  \Delta(|Ax_0|)=\mathop{\textup{argmin}}_{|Ax|=|Ax_0|}\sigma_k(x)_X
  $$
is one such decoder.
\eproof

\subsection{The Case $ X=\ell_1 $} We next apply Theorem \ref{lemmaSNSP} to the case $X=\ell_1$ norm.
The following lemma establishes a relation between S-RIP and S-NSP for the $ \ell_1 $-norm.

\begin{lemma} \label{relationl1}
 Let $ a, b, k $ be integers. Assume that $ A\in\mathbb{R}^{m\times N} $ satisfies the S-RIP of order $ (a+b)k $ with constants $ \theta_-, \ \theta_+\in(0, 2) $. Then $ A $ satisfies the S-NSP under the $\ell_1$-norm of order $ ak $ with constant
 $$
 C_0=1+\sqrt{\frac{a(1+\delta)}{b(1-\delta)}},
 $$
 where $ \delta $ is the restricted isometry constant and
 $\delta:=\max\{1-\theta_-,\theta_+-1\}<1$.
\end{lemma}

We remark that the above lemma is the analogous to the following lemma providing a relationship between RIP and NSP, which was shown in \cite{CDR09}:

\begin{lemma}[\cite{CDR09}, Lemma 4.1] \label{4.1in1}
  Let $a=l/k$, $b=l'/k$ where $l,l'\geq k$ are integers. Assume that $ A\in\mathbb{R}^{m\times N}$ satisfies the RIP of order $ (a+b)k $ with $ \delta=\delta_{(a+b)k}<1 $. Then $ A $ satisfies the null space property for the $ \ell_1 $-norm of order $ ak $ with constant $C_0=1+\frac{\sqrt{a(1+\delta)}}{\sqrt{b(1-\delta)}} $.
\end{lemma}

\vspace{3mm}

\noindent
{\bf Proof of Lemma \ref{relationl1}.}~~
By the definition of S-RIP, for any index set $ I\subseteq\{1,\ldots,m\} $ with $ |I|\geq m/2 $, the matrix $ A_I\in\R^{|I|\times N} $ satisfies the RIP of order $ (a+b)k $ with constant $\delta_{(a+b)k}=\delta:=\max\{1-\theta_-,\theta_+-1\}< 1 $. It follows from Lemma \ref{4.1in1} that
$$
\|\eta\|_1\leq\left(1+\sqrt{\frac{a(1+\delta)}{b(1-\delta)}}\,\right)\sigma_{ak}(\eta)_1
$$
for all $\eta\, \in\mathcal{N}(A_I)$. This proves the lemma.
\eproof

Set $ a=2$ and $b=1 $ in Lemma \ref{relationl1} we infer that for the $\ell_1$-norm if $ A $ satisfies the S-RIP of order $ 3k $ with constants $ \theta_-, \ \theta_+\in (0, 2)$, then $ A $ satisfies the S-NSP of order $ 2k $ with constant $C_0=1+\sqrt{\frac{2(1+\delta)}{1-\delta}} $. Hence by Theorem \ref{lemmaSNSP}, there must exist a decoder that has the instance optimality under the $ \ell_1 $-norm with constant $ 2C_0 $. According to Theorem \ref{gaussiansrip}, by taking $m=O(k\log (N/k))$ a Gaussian random matrix $A$ satisfies S-RIP of order $3k$ with high probability. Hence, there exists a decoder $\Delta$ so that the pair $(A, \Delta)$ has the the $\ell_1$-norm phaseless instance optimality at the cost of $m=O(k\log (N/k))$ measurements, as with the traditional instance optimality.

We are now ready to prove the following theorem on phaseless instance optimality under the $ X=\ell_1$-norm.

\begin{theorem}\label{theoreml1min}
 Let $ A \in\R^{m\times N}$ satisfy the S-RIP of order $ t\cdot k $  with constants $ 0<\theta_-<1< \theta_+<2 $, where
$$
    t\geq\max\left\lbrace \frac{2}{\theta_-}, \frac{2}{2-\theta_+} \right\rbrace>2.
$$
Let
\begin{equation}\label{l1minimation}
 \Delta(|Ax_0|)=\mathop{\textup{argmin}}_{x\in\R^N}\left\lbrace \|x\|_1: |Ax|=|Ax_0|\right\rbrace.
\end{equation}
 Then $ (A,\Delta) $ has the $\ell_1$-norm phaseless instance optimality with constant $ C=\frac{2C_0}{2-C_0} $, where $C_0=1+\sqrt{\frac{1+\delta}{(t-1)(1-\delta)}} $ and as before
$$
    \delta:=\max\{1-\theta_-,\theta_+-1\}\leq 1-\frac{2}{t}.
$$
\end{theorem}
\Proof
 Let $x_0 \in\R^N$ and set $\hat x =\Delta(|Ax_0|)$. Then by definition
$$
    \|\hat{x}\|_1\leq\|x_0\|_1\quad \text{and}\quad |A\hat{x}|=|Ax_0|.
$$
Denote by $I\subseteq \{1,\ldots,m\}$ the set of indices
$$
   I=\left\{j: \innerp{a_j,\hat{x}}=\innerp{a_j,x_0}\right\},
$$
and thus $\innerp{a_j,\hat{x}}=-\innerp{a_j,x_0}$ for $j \in I^c$. It follows that
$$
   A_I(\hat{x}-x_0)=0 \quad \mbox{and} \quad A_{I^c}(\hat{x}+x_0)=0.
$$
Set
$$
\eta:=\hat{x}-x_0\in\mathcal{N}(A_I).
$$
We know that $ A $ satisfies the S-RIP of order $ tk $ with constants
$ \theta_-,\ \theta_+ $ where
$$
   t\geq\max\left\lbrace \frac{2}{\theta_-}, \frac{2}{2-\theta_+} \right\rbrace>2.
$$

For the case $\abs{I}\geq m/2$, $A_I$ satisfies the RIP of order $tk$ with RIP constant
$$
    \delta=\delta_{tk}:=\max\{1-\theta_-, \theta_+-1\}\leq 1-\frac{2}{t}< 1.
$$
Take $ a:=1,\ b:=t-1 $ in Lemma \ref{relationl1}. Then $ A $ satisfies the $\ell_1$-norm S-NSP of order $ k $ with constant
$$
    C_0=1+\sqrt{\frac{1+\delta}{(t-1)(1-\delta)}}<2.
$$
This yields
\begin{equation}\label{theorem3.3-1}
   \|\eta\|_1\leq C_0\|\eta_{T^c}\|_1,
\end{equation}
where $ T $ is the index set for the $k$ largest coefficients of $ x_0 $ in magnitude.
Hence $ \|\eta_T\|_1\leq(C_0-1)\|\eta_{T^c}\|_1 $. Since $ \|\hat{x}\|_1\leq\|x_0\|_1 $ we have
\begin{align*}
   \|x_0\|_1\geq\|\hat{x}\|_1 &=\|x_0+\eta\|_1 =\|x_{0,T}+x_{0,T^c}+\eta_T+\eta_{T^c}\|_1\\
    &\geq\|x_{0,T}\|_1-\|x_{0,T^c}\|_1+\|\eta_{T^c}\|_1-\|\eta_T\|_1.
\end{align*}
It follows that
$$
    \|\eta_{T^c}\|_1\leq\|\eta_T\|_1+2\sigma_k(x_0)_1\leq(C_0-1)\|\eta_{T^c}\|_1+2\sigma_k(x_0)_1
$$
and thus
$$
     \|\eta_{T^c}\|_1\leq\frac{2}{2-C_0}\sigma_k(x_0)_1.
$$
Now (\ref{theorem3.3-1}) yields
$$
   \|\eta\|_1\leq C_0\|\eta_{T^c}\|_1\leq\frac{2C_0}{2-C_0}\sigma_k(x_0)_1,
$$
which implies
$$
  \|\hat{x}-x_0\|_1\leq C_0\|\eta_{T^c}\|_1\leq\frac{2C_0}{2-C_0}\sigma_k(x_0)_1.
$$

For the case $\abs{I^c}\geq m/2$ identical argument yields
 $$
 \|\hat{x}+x_0\|_1\leq C_0\|\eta_{T^c}\|_1\leq\frac{2C_0}{2-C_0}\sigma_k(x_0)_1.
 $$
The theorem is now proved.
\eproof

From Theorem \ref{gaussiansrip}, we know that an $m\times N$ random Gaussian matrix with $m=\mathcal{O}(tk\log(N/k))$ satisfies the S-RIP of order $ tk $ with high probability. We therefore conclude that the $\ell_1$-norm phaseless instance optimality of order $ k $ can be achieved at the cost of $m=\mathcal{O}(tk\log(N/k))$ measurements.

\subsection{Mixed-Norm phaseless Instance Optimality}

We now consider {\em mixed-norm  phaseless instance optimality}. Let $ 1\leq q\leq p\leq 2 $ and $s=1/q-1/p $. We seek estimates of the form
\begin{equation}\label{eq:pqinst2}
    \min\{\|x-\Delta(\abs{Ax})\|_{p}, \|x+\Delta(\abs{Ax})\|_{p}\}
        \leq C_0 k^{-s}\sigma_k(x)_{q}
\end{equation}
for all $x\in\R^N$. We shall prove both necessary and sufficient conditions for mixed-norm phaseless instance optimality.

\begin{theorem}\label{mixsuf}
    Let $A\in \R^{ m\times N} $ and $1\leq q\leq p\leq 2$. Set $s=1/q-1/p $. Then a decoder $ \Delta $ satisfying the mixed norm phaseless instance optimality (\ref{eq:pqinst2}) with constant $ C_0 $ exists if: for any index set $ I\subseteq\{1,\ldots,m\} $ and any $\eta_1\in\mathcal{N}(A_I)$, $\eta_2\in\mathcal{N}(A_{I^c})$ we have
\begin{equation}\label{mixeq:suf}
   \min\{\|\eta_1\|_p  ,  \|\eta_2\|_p \}
        \leq \frac{C_0}{4}k^{-s}\Bigl(\sigma_k(\eta_1-\eta_2)_q+\sigma_k(\eta_1+\eta_2)_q\Bigr).
\end{equation}
Conversely, assume a decoder $\Delta$ satisfying the mixed norm phaseless instance optimality (\ref{eq:pqinst2}) exists. Then for any index set $ I\subseteq\{1,\ldots,m\} $ and any $\eta_1\in\mathcal{N}(A_I)$, $\eta_2\in\mathcal{N}(A_{I^c})$ we have
\begin{equation*}
    \min\{\|\eta_1\|_p  ,  \|\eta_2\|_p \}
        \leq \frac{C_0}{2}k^{-s}\Bigl(\sigma_k(\eta_1-\eta_2)_q+\sigma_k(\eta_1+\eta_2)_q\Bigr).
\end{equation*}
\end{theorem}
\Proof   The proof is virtually identical to the proof of Theorem \ref{nsp_general}. We shall omit the details here in the interest of brevity.
\eproof

\begin{definition}(Mixed-Norm Strong Null Space Property)
{\rm We say that $ A $ has the mixed strong null space property in norms $ (\ell_p,\ell_q) $ of order $ k $ with constant $ C $ if for any index set $ I\subseteq\{1,\ldots,m\} $ with $ |I|\geq m/2 $ the matrix $ A_I\in\R^{|I|\times N} $ satisfies
$$
    \|\eta\|_p\leq Ck^{-s}\sigma_k(\eta)_q
$$
for all $\eta\in\mathcal{N}(A_I)$, where $s = 1/q-1/p$.
}
\end{definition}

The above is an extension of the standard definition of the mixed null space property of order $k$ in norms $ (\ell_p,\ell_q) $ (see \cite{CDR09}) for a matrix $A$, which requires
$$
    \|\eta\|_p\leq Ck^{-s}\sigma_k(\eta)_q
$$
for all $\eta\in\mathcal{N}(A)$. We have the following straightforward generalization of Theorem \ref{lemmaSNSP}.

\begin{theorem}\label{mixSNSP}
   Assume that $ A\in\mathbb{R}^{m\times N} $ has the mixed strong null space property of order $2k$ in norms $ (\ell_p,\ell_q) $ with constant $ C_0/2 $, where $ 1\leq q\leq p\leq 2 $. Then there exists a decoder $\Delta$ such that the mixed-norm phaseless instance optimality (\ref{eq:pqinst2}) holds with constant $C_0 $.
\end{theorem}

We establish relationships between mixed-norm strong null space property and the S-RIP. First we present the following lemma that was proved in \cite{CDR09}.

\begin{lemma}[\cite{CDR09}, Lemma 8.2]  \label{8.2 in 1}
    Let $k\geq 1$ and $\tilde k = k(\frac{N}{k})^{2-2/q}$. Assume that $A\in\R^{m\times N}$ satisfies the RIP of order $ 2k+\tilde{k}$ with $ \delta:=\delta_{2k+\tilde{k}}<1 $. Then $A$ satisfies the mixed null space property in norms $(\ell_p,\ell_q)$ of order $ 2k $ with constant $ C_0=2^{1/p+1/2}\sqrt{\frac{1+\delta}{1-\delta}}+2^{1/p-1/q}$.
\end{lemma}

\vspace{2mm}

\begin{prop}\label{mixsnsprela}
 Let $k\geq 1$ and $\tilde k = k(\frac{N}{k})^{2-2/q}$. Assume that $A\in\R^{m\times N}$ satisfies the S-RIP of order $ 2k+\tilde{k}$ with constants  $ 0<\theta_- <1 <\theta_+<2$. Then $ A $ satisfies the mixed strong null space property in norms $ (\ell_p, \ell_q) $ of order $ 2k $ with constant $C_0=2^{1/p+1/2}\sqrt{\frac{1+\delta}{1-\delta}}+2^{1/p-1/q} $, where $ \delta $ is  the RIP constant and $\delta:=\delta_{2k+\tilde{k}}= \max\{1-\theta_-, \theta_+-1\}$.
\end{prop}
\Proof
By definition for any index set $ I\subseteq\{1,\ldots,m\} $ with $ |I|\geq m/2 $, the matrix $ A_I\in\R^{|I|\times N} $ satisfies RIP of order $ 2k+\tilde{k} $ with constant $ C_0=2^{1/p+1/2}\sqrt{\frac{1+\delta}{1-\delta}}+2^{1/p-1/q} $,  where $ \delta $ is  the RIP constant and $\delta:=\delta_{2k+\tilde{k}}= \max\{1-\theta_-, \theta_+-1\}$. By Lemma \ref{8.2 in 1}, we know that $ A_I $ satisfies the mixed null space property in norms $ (\ell_p,\ell_q) $ of order $ 2k $ with constant $ C_0=2^{1/p+1/2}\sqrt{\frac{1+\delta}{1-\delta}}+2^{1/p-1/q} $, in other words for any $\eta\in\mathcal{N}(A_I)$,
$$
    \|\eta\|_p\leq Ck^{-s}\sigma_{2k}(\eta)_q.
$$
So $ A $ satisfies the mixed strong null space property.
\eproof

\vspace{2mm}

\begin{corollary} \label{co:41}
Let $k\geq 1$ and $\tilde k = k(\frac{N}{k})^{2-2/q}$. Assume that $A\in\R^{m\times N}$ satisfies the S-RIP of order $ 2k+\tilde{k}$ with constants  $ 0<\theta_- <1 <\theta_+<2$. Let $\delta:=\delta_{2k+\tilde{k}} =\max\{1-\theta_-, \theta_+-1\}<1 $. Define the decoder $ \Delta $  for $ A $ by
\begin{equation}
    \Delta(|Ax_0|)=\mathop{\textup{argmin}}_{|Ax|=|Ax_0|}\sigma_k(x)_q.
\end{equation}
Then (\ref{eq:pqinst2}) holds with constant $2C_0$, where
$C_0=2^{1/p+1/2}\sqrt{\frac{1+\delta}{1-\delta}}+2^{1/p-1/q}  $.
\end{corollary}
\Proof
    By the Proposition \ref{mixsnsprela}, the matrix $ A $ satisfies the mixed strong null space property in $ (\ell_p,\ell_q) $ of order $ 2k $ with constant $ C_0=2^{1/p+1/2}\sqrt{\frac{1+\delta}{1-\delta}}+2^{1/p-1/q} $. The corollary now follows immediately from  Theorem \ref{mixSNSP}.
\eproof

\begin{remark}{\rm
    Combining Theorem \ref{gaussiansrip} and Corollary \ref{co:41}, the mixed phaseless  instance optimality of order $ k $ in norms $ (\ell_p,\ell_q) $ can be achieved at the price of $ \mathcal{O}(k(N/k)^{2-2/q}\log(N/k)) $ measurements, just as with the traditional mixed $(\ell_p,\ell_q)$-norm instance optimality. Theorem \ref{maintherom} implies that the $\ell_1$ decoder satisfies the $(p,q)=(2,1)$ mixed-norm phaseless instance optimality  at the price of $ \mathcal{O}(k\log(N/k)) $ measurements.
}
\end{remark}

\section{Appendix: Proof of Lemma \ref{mainlemma}}

\setcounter{equation}{0}

We will first need the following two Lemmas to prove Lemma \ref{mainlemma}.

\begin{lemma}[Sparse Representation of a Polytope \cite{CZ14,XX13}]  \label{quotelemma1}
Let $s\geq 1$ and  $\alpha>0$. Set
$$
   T(\alpha,s):=\Bigl\{u\in\mathbb{R}^n: \|u\|_\infty\leq\alpha,\ \|u\|_1\leq s\alpha\Bigr\}.
$$
For any $v\in\mathbb{R}^n$ let
$$
  \mathit{U}(\alpha,s,v):=\Bigl\{u\in\mathbb{R}^n:\textup{supp}(u)\subseteq
          \textup{supp}(v),\|u\|_0\leq s,\|u\|_1=\|v\|_1,\|u\|_\infty\leq\alpha\Bigr\}.
$$
Then $v\in T(\alpha,s)$ if and only if $v$ is in the convex hull of $\mathit{U}(\alpha,s,v)$, i.e. $v$ can be expressed as a convex combination of some $u_1, \dots, u_N$ in $\mathit{U}(\alpha,s,v)$.
\end{lemma}

\begin{lemma}[Lemma 5.3 in \cite{ref6}]    \label{quotelemma2}
  Assume that $ a_1\geq a_2\geq\cdots\geq a_m\geq 0 $. Let $r \leq m$ and $\lambda \geq 0$ such that $ \sum_{i=1}^{r}a_i + \lambda\geq\sum_{i=r+1}^{m}a_i $. Then for all $ \alpha\geq 1 $ we have
\begin{equation}  \label{l1l2relation}
   \sum_{j=r+1}^{m}a_j^\alpha\leq r
         \left(\sqrt[\alpha]{\frac{\sum_{i=1}^{r}a_i^\alpha}{r}}+\frac{\lambda}{r}\right)^\alpha.
\end{equation}
In particular for $\lambda=0$ we have
$$
     \sum_{j=r+1}^{m}a_j^\alpha\leq\sum_{i=1}^{r}a_i^\alpha.
$$
\end{lemma}

We are now ready to prove Lemma \ref{mainlemma}.

\vspace{2mm}
\noindent
{\bf Proof of Lemma \ref{mainlemma}.}~~ Set $h:=\hat{x}-x_0$. Let $T_0$ denote the set of the largest $k$ coefficients of $x_0$ in magnitude. Then
\begin{align*}
      \|x_0\|_1+\rho & \geq\|\hat{x}\|_1 =\|x_0+h\|_1 \\
                     &=\|x_{0,T_0}+h_{T_0}+x_{0,T_0^c}+h_{T_0^c}\|_1\\
                     &\geq \|x_{0,T_0}\|_1-\|h_{T_0}\|_1-\|x_{0,T_0^c}\|_1+\|h_{T_0^c}\|_1.
\end{align*}
It follows that
\begin{align*}
     \|h_{T_0^c}\|_1&\leq\|h_{T_0}\|_1+2\|x_{0,T_0^c}\|_1+\rho\\
      &=\|h_{T_0}\|_1+2\sigma_k(x_0)_1+\rho.
 \end{align*}
   Suppose that $ S_0 $ is the index set of the $k$ largest entries in absolute value of $ h $.  Then we can get
\begin{align*}
    \|h_{S_0^c}\|_1\leq\|h_{T_0^c}\|_1&\leq\|h_{T_0}\|_1+2\sigma_k(x_0)_1+\rho
       \\ &\leq\|h_{S_0}\|_1+2\sigma_k(x_0)_1+\rho.
\end{align*}
Set
$$
    \alpha:=\frac{\|h_{S_0}\|_1+2\sigma_k(x_0)_1+\rho}{k}.
$$
We divide $ h_{S_0^c }$ into two parts $ h_{S_0^c }=h^{(1)}+h^{(2)} $, where
\begin{align*}
    h^{(1)}:=h_{S_0^c}\cdot I_{\{i:\,|h_{S_0^c}(i)|>\alpha/(t-1)\}}, \quad
    h^{(2)}:=h_{S_0^c}\cdot I_{\{i:\,|h_{S_0^c}(i)|\leq\alpha/(t-1)\}} .
\end{align*}
 A simple observation is that $\|h^{(1)}\|_1\leq\|h_{S_0^c}\|_1\leq\alpha k $.   Set
$$
  \ell := |\supp(h^{(1)})|=\|h^{(1)}\|_0 .
$$
Since all non-zero entries of $ h^{(1)} $ have magnitude larger than $ \alpha/(t-1) $, we have
$$
   \alpha k\geq \|h^{(1)}\|_1=\sum_{i\in \textup{supp}(h^{(1)})}|h^{(1)}(i)|
      \geq\sum_{i\in \textup{supp}(h^{(1)})}\frac{\alpha}{t-1}=\frac{\alpha \ell}{t-1},
$$
which implies  $ \ell\leq (t-1)k $.  Thus we have:
\begin{equation}\label{opt-rip}
    \big\langle A(h_{S_0}+h^{(1)}), Ah\big\rangle \leq\|A(h_{S_0}+h^{(1)})\|_2\cdot \|Ah\|_2
    \leq\sqrt{1+\delta}\cdot \|h_{S_0}+h^{(1)}\|_2\cdot\epsilon.
\end{equation}
Here we apply the facts that $\|h_{S_0}+h^{(1)}\|_0=\ell+k\leq tk $ and  $ A $ satisfies the RIP of order $ tk $ with $ \delta:=\delta_{tk}^A $. We shall assume at first that $ tk $ as an integer. Note that $\|h^{(2)}\|_\infty \leq \frac{\alpha}{t-1}$ and
\begin{equation}
    \|h^{(2)}\|_1 =\|h_{S_0^c}\|_1-\|h^{(1)}\|_1
         \leq k\alpha- \frac{ \alpha\ell}{t-1}=(k(t-1)-\ell)\frac{\alpha}{t-1}.
\end{equation}
We take $ s:=k(t-1)-\ell $ in Lemma \ref{quotelemma1}  and obtain  that  $ h^{(2)} $ is a weighted mean
$$
    h^{(2)}=\sum_{i=1}^{N}\lambda_iu_i,\quad \quad 0\leq \lambda_i\leq 1, \quad \sum_{i=1}^N\lambda_i=1
$$
where $ \|u_i\|_0\leq k(t-1)-\ell,  \|u_i\|_1=\|h^{(2)}\|_1 $, $\|u_i\|_\infty\leq\alpha/(t-1) $ and $\textup{supp}(u_i)\subseteq \textup{supp}(h^{(2)}) $. Hence
\begin{align*}
       \|u_i\|_2\leq\sqrt{\|u_i\|_0}\cdot \|u_i\|_\infty &
                 =\sqrt{k(t-1)-\ell }\cdot \|u_i\|_\infty\\
                 & \leq\sqrt{k(t-1)}\cdot \|u_i\|_\infty\\
                 &\leq\alpha\sqrt{k/(t-1)}.
\end{align*}
Now for $ 0\leq\mu\leq 1 $ and $ d\geq 0$, which will be chosen later, set
$$
      \beta_j:=h_{S_0}+h^{(1)}+\mu\cdot u_j, \quad j=1,\ldots,N.
$$
Then for fixed $i\in[1,N]$
\begin{align*}
      \sum_{j=1}^{N}\lambda_j\beta_j-d\beta_i&=h_{S_0}+h^{(1)}+\mu\cdot h^{(2)}-d\beta_i\\
          &=(1-\mu-d)(h_{S_0}+h^{(1)})-d\mu u_i+\mu h.
\end{align*}
Recall that  $\alpha=\frac{\|h_{S_0}\|_1+2\sigma_k(x_0)_1+\rho}{k}$.  Thus
\begin{align} \label{opt-u}
   \|u_i\|_2&\leq\sqrt{k/(t-1)}\alpha \\\nonumber
    &\leq\frac{\|h_{S_0}\|_2}{\sqrt{t-1}}+\frac{2\sigma_k(x_0)_1+\rho}{\sqrt{k(t-1)}}\\\nonumber
     &\leq\frac{\|h_{S_0}+h^{(1)}\|_2}{\sqrt{t-1}}+\frac{2\sigma_k(x_0)_1+\rho}{\sqrt{k(t-1)}}\\\nonumber
      &=\frac{z+R}{\sqrt{t-1}},\nonumber
\end{align}
where   $ z:=\|h_{S_0}+h^{(1)}\|_2$ and $R:=\frac{2\sigma_k(x_0)_1+\rho}{\sqrt{k}}$. It is easy to check the following identity:
\begin{align}\label{eq:lamdabeta2}
  (2d-1)&\sum_{1\leq i<j\leq N}\lambda_i\lambda_j\|A(\beta_i-\beta_j)\|_2^2 \nonumber \\
        &=\sum_{i=1}^{N}\lambda_i\Bigl\|A(\sum_{j=1}^{N}\lambda_j\beta_j-d\beta_i)\Bigr\|_2^2
        - \sum_{i=1}^{N}\lambda_i(1-d)^2\|A\beta_i\|_2^2,
\end{align}
provided that $\sum_{i=1}^N\lambda_i=1$. Choose  $d=1/2$
in (\ref{eq:lamdabeta2}) we then have
$$
    \sum_{i=1}^{N}\lambda_i\Bigl\|A\Bigl( (\frac{1}{2}-\mu)(h_{S_0}+h^{(1)})-\frac{\mu} {2}u_i+\mu h\Bigr)\Bigr\|_2^2-\sum_{i=1}^{N}\frac{\lambda_i}{4}\|A\beta_i\|_2^2 = 0.
$$
Note that for $d=1/2$,
\begin{align*}
  \Bigl\|A\Bigl( &(\frac{1}{2}-\mu)(h_{S_0}+h^{(1)})-\frac{\mu} {2}u_i+\mu h\Bigr)\Bigr\|_2^2\\
  &= \Bigl\|A\Bigl( (\frac{1}{2}-\mu)(h_{S_0}+h^{(1)})-\frac{\mu}{2}u_i\Bigr) \Bigr\|_2^2 +2\Bigl\langle A\Bigl((\frac{1}{2}-\mu)(h_{S_0}+h^{(1)})-\frac{\mu}{2}u_i\Bigr), \mu Ah\Bigr\rangle +\mu^2\|Ah\|_2^2.
\end{align*}
It follows from $\sum_{i=1}^N \lambda_i =1$ and $h^{(2)}=\sum_{i=1}^{N}\lambda_iu_i$ that
\begin{align}
 &\quad \sum_{i=1}^{N}\lambda_i\Bigl\|A\Bigl( (\frac{1}{2}-\mu)(h_{S_0}+h^{(1)})-\frac{\mu} {2}u_i+\mu h\Bigr)\Bigr\|_2^2\nonumber \\
 &= \sum_i\lambda_i\Bigl\|A\Bigl( (\frac{1}{2}-\mu)(h_{S_0}+h^{(1)})-\frac{\mu}{2}u_i\Bigr) \Bigr\|_2^2 +2\Bigl\langle A\Bigl((\frac{1}{2}-\mu)(h_{S_0}+h^{(1)})-\frac{\mu}{2}h^{(2)}\Bigr), \mu Ah\Bigr\rangle +\mu^2\|Ah\|_2^2 \nonumber\\
 & = \sum_i\lambda_i\Bigl\|A\Bigl( (\frac{1}{2}-\mu)(h_{S_0}+h^{(1)})-\frac{\mu}{2}u_i\Bigr) \Bigr\|_2^2
     + \mu(1-\mu)\Bigl\langle
    A(h_{S_0}+h^{(1)}),Ah\Bigr\rangle -\sum_{i=1}^{N}\frac{\lambda_i}{4}\|A\beta_i\|_2^2.  \label{eq:threeterm}
\end{align}
Set $ \mu=\sqrt{t(t-1)}-(t-1)$. We next estimate the three terms in (\ref{eq:threeterm}). Noting that $\| h_{S_0}\|_0\leq k $, $ \|h^{(1)}\|_0\leq \ell$ and  $\|u_i\|_0\leq s =k(t-1)-\ell$, we obtain
$$
      \| \beta_i\|_0\leq  \| h_{S_0}\|_0 +  \|h^{(1)}\|_0+ \|u_i\|_0\leq t\cdot k
$$
and  $\| (\frac{1}{2}-\mu)(h_{S_0}+h^{(1)})-\frac{\mu}{2}u_i\|_0\leq t\cdot k $. Since $A$ satisfies the RIP of order $t\cdot k$ with $\delta $, we have
\begin{align*}
   \Bigl\|A\Bigl( (\frac{1}{2}-\mu)(h_{S_0}+h^{(1)})-\frac{\mu}{2}u_i\Bigr) \Bigr\|_2^2&\leq (1+\delta)\|(\frac{1}{2}-\mu)(h_{S_0}+h^{(1)})-\frac{\mu}{2}u_i\|_2^2\\
       &= (1+\delta)\Bigl((\frac{1}{2}-\mu)^2 \|(h_{S_0}+h^{(1)})\|_2^2+\frac{\mu^2}{4}\|u_i\|_2^2\Bigr)\\
       &= (1+\delta)\Bigl((\frac{1}{2}-\mu)^2 z^2+\frac{\mu^2}{4}\|u_i\|_2^2\Bigr)
\end{align*}
and
$$
   \|A\beta_i\|_2^2\geq (1-\delta)\|\beta_i\|_2^2 =(1-\delta)( \|h_{S_0}+h^{(1)}\|_2^2+\mu^2\cdot \|u_i\|_2^2)=
   (1-\delta)(z^2+\mu^2\cdot \|u_i\|_2^2).
$$
Combining the result above with (\ref{opt-rip}) and (\ref{opt-u}) we  get
\begin{align*}
    0&\leq(1+\delta)\sum_{i=1}^{N}\lambda_i\Bigl((\frac{1}{2}-\mu)^2z^2+\frac{\mu^2}{4}\|u_i\|_2^2\Bigr)
    +\mu(1-\mu)\sqrt{1+\delta}\cdot z\cdot\epsilon-(1-\delta)\sum_{i=1}^{N}\frac{\lambda_i}{4}(z^2+\mu^2\|u_i\|_2^2)\\
    &=\sum_{i=1}^{N}\lambda_i\Bigl(\Bigl((1+\delta)(\frac{1}{2}-\mu)^2-\frac{1-\delta}{4} \Bigr)z^2+\frac{\delta}{2}\mu^2\|u_i\|_2^2  \Bigr)+\mu(1-\mu)\sqrt{1+\delta}\cdot z\cdot\epsilon \\
    &\leq\sum_{i=1}^{N}\lambda_i\Bigl(\Bigl((1+\delta)(\frac{1}{2}-\mu)^2-\frac{1-\delta}{4} \Bigr)z^2+\frac{\delta}{2}\mu^2\frac{(z+R)^2}{t-1}  \Bigr)+\mu(1-\mu)\sqrt{1+\delta}\cdot z\cdot\epsilon \\
    &=\Bigl((\mu^2-\mu)+\delta\Bigl( \frac{1}{2}-\mu+(1+\frac{1}{2(t-1)})\mu^2\Bigr)  \Bigr)z^2+\Bigl( \mu(1-\mu)\sqrt{1+\delta}\cdot\epsilon+\frac{\delta\mu^2R}{t-1}\Bigr)z+\frac{\delta\mu^2R^2}{2(t-1)} \\
    &=-t\Bigl((2t-1)-2\sqrt{t(t-1)} \Bigr) (\sqrt{\frac{t-1}{t}}-\delta)z^2+\Bigl( \mu^2\sqrt{\frac{t}{t-1}}\sqrt{1+\delta}\cdot\epsilon+\frac{\delta\mu^2R}{t-1}\Bigr)z+\frac{\delta\mu^2R^2}{2(t-1)}\\
    &=\frac{\mu^2}{t-1}\Bigl(-t(\sqrt{\frac{t-1}{t}}-\delta)z^2+(\sqrt{t(t-1)(1+\delta)}\epsilon+\delta R)z+\frac{\delta R^2}{2} \Bigr),
\end{align*}
which is a quadratic inequality for $z$. We know $\delta<\sqrt{(t-1)/t}$. So by solving the above inequality we get
\begin{align*}
    z&\leq\frac{(\sqrt{t(t-1)(1+\delta)}\epsilon+\delta R)+\left((\sqrt{t(t-1)(1+\delta)}\epsilon+\delta R)^2+2t(\sqrt{(t-1)/t}-\delta)\delta R^2 \right)^{1/2}  }
    {2t(\sqrt{(t-1/t)}-\delta) }\\
    &\leq\frac{\sqrt{t(t-1)(1+\delta)}}{t(\sqrt{(t-1)/t}-\delta)}\epsilon+\frac{2\delta+\sqrt{2t(\sqrt{(t-1)/t}-\delta)\delta}}{2t(\sqrt{(t-1)/t}-\delta)}R.
\end{align*}
Finally, noting that  $ \|h_{S_0^c}\|_1\leq\|h_{S_0}\|_1+R\sqrt{k} $, in the Lemma \ref{quotelemma2}, if we set $ m=N $, $ r=k $, $ \lambda=R\sqrt{k}\geq 0 $ and $ \alpha=2 $ then $\|h_{S_0^c}\|_2\leq\|h_{S_0}\|_2+R $. Hence
\begin{align*}
    \|h\|_2&=\sqrt{\|h_{S_0}\|_2^2+\|h_{S_0^c}\|_2^2}\\
    &\leq\sqrt{\|h_{S_0}\|_2^2+(\|h_{S_0}\|_2+R)^2}\\
    &\leq\sqrt{2\|h_{S_0}\|_2^2}+R\leq\sqrt{2}z+R\\
    &\leq\frac{\sqrt{2(1+\delta)}}{1-\sqrt{t/(t-1)}\delta}\epsilon+\left( \frac{\sqrt{2}\delta+\sqrt{t(\sqrt{(t-1)/t}-\delta)\delta}}{t(\sqrt{(t-1)/t}-\delta)}+1 \right) R.
\end{align*}
Substitute $R$ into this inequality and the conclusion follows.

     For the case where $t\cdot k$ is not an integer, we set $t^*:=\lceil tk\rceil / k$, then $t^*>t$ and $\delta_{t^*k}=\delta_{tk}<\sqrt{\frac{t-1}{t}}<\sqrt{\frac{t^*-1}{t^*}}$. We can then prove the result by working on $\delta_{t^*k}$.
\eproof

 \end{document}